\documentclass{agtart_a}
\pdfoutput=1
\usepackage{pinlabel}

\title[Legendrian links and the spanning tree model for Khovanov homology]{Legendrian links and the spanning tree model\\for Khovanov homology}

\author{Hao Wu}
\givenname{Hao}
\surname{Wu}
\address{Department of mathematics and Statistics\\
Lederle Graduate Research Tower\\\newline
710 North Pleasant Street\\
University of Massachusetts\\Amherst, MA 01003-9305\\USA}
\email{wu@math.umass.edu}
\urladdr{}

\volumenumber{6}
\issuenumber{}
\publicationyear{2006}
\papernumber{61}
\startpage{1745}
\endpage{1757}

\doi{}
\MR{}
\Zbl{}

\keyword{Legendrian link}
\keyword{Khovanov homology}
\keyword{Tait graph}
\keyword{spanning tree}
\subject{primary}{msc2000}{57M25}
\subject{secondary}{msc2000}{57R17}

\received{31 May 2006}
\revised{25 August 2006}
\accepted{26 August 2006}
\published{21 October 2006}
\publishedonline{21 October 2006}
\proposed{}
\seconded{}
\corresponding{}
\editor{CPR}
\version{}

\arxivreference{math.GT/0605630}

\AtBeginDocument{\let\bar\wbar\let\tilde\wtilde\let\hat\what}
\numberwithin{equation}{section}
\makeop{cut}
\makeop{cyc}

\newcommand{\zed}{\mathbb{Z}}
\newcommand{\Real}{\mathbb{R}}

\makeatletter
\def\cnewtheorem#1[#2]#3{\newtheorem{#1}{#3}[section]
\expandafter\let\csname c@#1\endcsname\c@theorem}
\makeatother

\theoremstyle{plain}
\newtheorem{theorem}{Theorem}
\numberwithin{theorem}{section}
\cnewtheorem{lemma}[theorem]{Lemma}
\cnewtheorem{proposition}[theorem]{Proposition}
\cnewtheorem{corollary}[theorem]{Corollary}
\cnewtheorem{conjecture}[theorem]{Conjecture}

\theoremstyle{definition}
\cnewtheorem{definition}[theorem]{Definition}
\cnewtheorem{example}[theorem]{Example}
\cnewtheorem{notation}[theorem]{Notation}
\cnewtheorem{observation}[theorem]{Observation}
\cnewtheorem{remark}[theorem]{Remark}
\cnewtheorem{question}[theorem]{Question}
\newtheorem*{acknowledgments}{Acknowledgments}

\newcommand{\rank}{\mathrm{rank}}

\begin{document}

\begin{abstract}
We use the spanning tree model for Khovanov homology to study
Legendrian links.  This leads to an alternative proof for Ng's
Khovanov bound for the Thurston--Bennequin number and to both a
necessary and a sufficient condition for this bound to be sharp.
\end{abstract}
\maketitle

\section{Introduction}\label{intro}

In \cite{K1}, M. Khovanov constructed a categorification of the Jones polynomial. That is, to any oriented link $L$, he associated a bigraded homology group $\mathcal{H}(L)$, the Khovanov homology, whose graded Euler characteristic is
\[
\chi(\mathcal{H}(L)) := \sum_{i,j}(-1)^iq^j\rank(\mathcal{H}^{i,j}(L)) = (q+q^{-1})V_L(q^2),
\]
where $V_L$ is the Jones polynomial, $i$ is the homological grading of $\mathcal{H}(L)$, and $j$ is the quantum grading of $\mathcal{H}(L)$.

The Khovanov homology has led to many interesting new developments in
knot theory and related fields. See Lee \cite{Lee1,Lee2}, Ng
\cite{NgKB}, Plamenevskaya \cite{Pl4} and Rasmussen \cite{Ras} for
examples. It is still very difficult to compute the Khovanov homology
in general. Recently, A Champanerkar and I Kofman \cite{CKVST} and,
independently, S Wehrli \cite{WehST} constructed a spanning tree
model for the Khovanov homology based on the spanning tree expansion
of the Jones polynomial introduced by M Thistlethwaite in
\cite{This}. Though the spanning tree model does not completely
determine the Khovanov homology, it does greatly simplify the Khovanov
chain complex used to compute the Khovanov homology. In some cases,
such simplifications are enough to deduce interesting results. For
example, Lee's result on the Khovanov homology of alternating knots is
reproved in \cite{CKVST,WehST} by the spanning tree model.

In this paper, we will use the spanning tree model for Khovanov homology to study Legendrian links in the standard contact $\Real^3$. In particular, we give an alternative proof of the following theorem of Ng.

\begin{theorem}[Ng \cite{NgKB}]\label{NgBound}
For any Legendrian link $L$ in the standard contact $\Real^3$,
\begin{equation}\label{inequality}
tb(L) \leq \min\{k~|\bigoplus_{j-i=k}\mathcal{H}^{i,j}(L)\neq 0\},
\end{equation}
where $tb$ is the Thurston--Bennequin number.
\end{theorem}

From our proof of \fullref{NgBound}, it's easy to see that we have the following necessary condition and sufficient condition for Ng's Khovanov bound to be sharp, where good spanning trees and bad spanning trees will be defined in \fullref{legendrian}.

\begin{theorem}\label{sharpness}
Let $L$ be a Legendrian link.

{\rm(i)}\qua If Inequality \eqref{inequality} is an equality, then the front projection of $L$ admits a good spanning tree.

{\rm(ii)}\qua If there is an integer $v$, such that the front projection of $L$ admits more good $v$--spanning trees than bad $(v+2)$--spanning trees, then Inequality \eqref{inequality} is an equality.
\end{theorem}

Specially, this theorem implies that Ng's bound is sharp for alternating links.

\begin{corollary}[Ng \cite{NgKB}]\label{alternating}
If $\mathcal{L}$ is an alternating link, then
\begin{equation}\label{alternating_sharp}
\overline{tb}(\mathcal{L}) = \min\{k~|\bigoplus_{j-i=k}\mathcal{H}^{i,j}(\mathcal{L})\neq0\},
\end{equation}
where $\overline{tb}(\mathcal{L})$ is the maximal Thurston--Bennequin number for a Legendrian link\break smoothly isotopic to $\mathcal{L}$.
\end{corollary}

It is interesting to compare \fullref{sharpness} to D Rutherford's results in \cite{Ruth}, where he demonstrated that the Kauffman polynomial bound for the Thurston--Bennequin number is sharp if and only if the front projection admits a ruling.

\begin{question}
Can we refine \fullref{sharpness} to get a necessary and sufficient condition for Ng's Khovanov bound to be sharp in terms of spanning trees?
\end{question}

\begin{acknowledgments}
The author would like to thank Abhijit Champanerkar, Ilya Kofman and Yongwu Rong for helpful discussions. He would also like to thank the referee for many helpful suggestions.
\end{acknowledgments}

\section{The spanning tree model for Khovanov homology}\label{model}

In this section, we recall the construction of the spanning tree model of Khovanov homology in \cite{CKVST}. For a similar construction, see also \cite{WehST}.

\begin{figure}[ht!]
\centering
\labellist\small
\pinlabel $+$ [l] at 145 41
\pinlabel $-$ [l] at 435 41
\endlabellist
\includegraphics[width=4in]{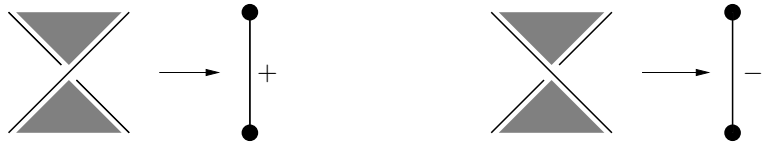}
\caption{Sign of an edge}\label{sign}
\end{figure}

Let $D$ be an oriented link diagram with an given ordering of
crossings. Checkerboard color the complementary regions of $D$. To
each black region, assign a vertex, and, to each crossing, assign an
edge connecting the two vertices corresponding to the two black
regions incident on this crossing. The result is a planar graph $G$
called a Tait graph of the link diagram. The edges of $G$ are
ordered by the ordering of the crossings in $D$. Assign a sign to
each edge of $G$ by the convention in \fullref{sign}.

In the rest of this section, we assume that the knot diagram $D$ is connected as a subset of $\mathbb{R}^2$. (This can be easily arranged using Reidemeister moves.)

Let $T$ be a spanning tree of $G$. For an edge $e\in T$, removing $e$ from $T$ divides $T$ into two connected components. Let $\cut(T,e)$ be the set of edges of $G$ connecting these two connected components of $T\setminus e$. $e$ is  said to be internally active if it has the lowest ordering among the elements of $\cut(T,e)$. A positive internally active edge is denoted by $L$, and a negative internally active edge is denoted by $\overline{L}$. An edge in $T$ that is not internally active is said to be internally inactive. A positive internally inactive edge is denoted by $D$, and a negative internally inactive edge is denoted by $\overline{D}$. For an edge $f\notin T$, $T\cup f$ contains a unique simple cycle. Let $\cyc(T,f)$ be the set of edges in this simple cycle. $f$ is said to be externally active if it has the lowest ordering among the elements of $\cyc(T,f)$. A positive externally active edge is denoted by $l$, and a negative externally active edge is denoted by $\overline{l}$. An edge outside $T$ that is not externally active is said to be externally inactive. A positive externally inactive edge is denoted by $d$, and a negative externally inactive edge is denoted by $\overline{d}$. Note that, for edges $e$ and $f$ with $e\in T$ and $f\notin T$, $f\in \cut(T,e)$ if and only if $e\in \cyc(T,f)$.

\begin{figure}[ht!]
\centering
\labellist
\pinlabel $A$ [b] at 113 38
\pinlabel $B$ [b] at 251 38
\endlabellist
\includegraphics[width=3in]{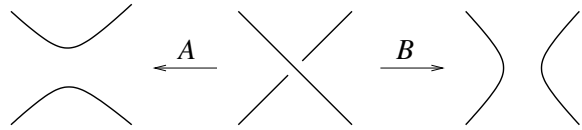}
\caption{Splicings of a crossing}\label{splicings}
\end{figure}

For any crossing in $D$, there are two ways to splice it, which are called the $A$--splicing and the $B$--splicing. These are depicted in \fullref{splicings}.

Given a spanning tree $T$ of the Tait graph $G$, one obtains a twisted unknot $U_T$ by splicing each inactive crossing following the rules in \fullref{splicingdead}. (cf \cite{CKVST}, Proposition 2.) Here, a twisted unknot is a diagram of the unknot obtained from the round circle by only Reidemeister I moves.

\begin{table}[ht!]
\begin{center}
\begin{tabular}{|c|c|c|c|}
\hline
\vrule width 0pt height 15pt
$D$ & $d$ & $\overline{D}$ & $\overline{d}$ \\
\hline
\vrule width 0pt height 15pt
$A$ & $B$ & $B$ & $A$ \\
\hline
\end{tabular}
\end{center}
\vspace{.4cm} \caption{Splicing inactive
crossings}\label{splicingdead}
\end{table}

The contribution of each active crossing to the writhe of $U_T$ is given in \fullref{UTwrithe}.

\begin{table}[ht!]
\begin{center}
\begin{tabular}{|c|c|c|c|}
\hline
\vrule width 0pt height 15pt
$L$ & $l$ & $\overline{L}$ & $\overline{l}$ \\
\hline
\vrule width 0pt height 15pt
$-$ & $+$ & $+$ & $-$ \\
\hline
\end{tabular}
\end{center}
\vspace{.4cm} \caption{Contribution of an active crossing to the
writhe of $U_T$}\label{UTwrithe}
\end{table}

Denote by $E_{\pm}(G)$ the number of positive/negative edges of $G$, by $V(G)$ the number of vertices of $G$. For type $X$ of edges, denote by $\#X_G(T)$ the number of edges of type $X$ in $G$ computed using $T$. When $G$ is clear from the context, we drop it from the notation.

For a spanning tree $T$ of $G$, define
\begin{eqnarray*}
u(T) & = & -w(U_T) = \#L(T)-\#l(T)-\#\overline{L}(T)+\#\overline{l}(T), \\
v(T) & = & E_+(T) + E_-(G\setminus T)= \#L(T)+\#D(T)+\#\overline{l}(T)+\#\overline{d}(T).
\end{eqnarray*}
where $w(\ast)$ means the writhe. (Note that the normalization of $v$ here is different from that in \cite{CKVST}.)

Define $\mathfrak{C}_T$ to be the bigraded free abelian group of rank two with one generator of bidegree $(u(T),v(T))$, and the other of bidegree $(u(T)+2,v(T)+2)$. In the rest of this paper, we call the first grading the $u$--grading, and the second the $v$--grading.

Define
\[
\mathfrak{C}(D)=\bigoplus_{T} \mathfrak{C}_T,
\]
where the direct sum is taken over all spanning trees of $G$.

\begin{theorem}[Champanerkar and Kofman \cite{CKVST}]\label{STmodel}
There is a boundary map $\partial$ of bidegree $(-1,-2)$ on $\mathfrak{C}(D)$, so that $(\mathfrak{C}(D),\partial)$ is a deformation retract of the Khovanov chain complex. In particular, we have
\[
\mathcal{H}^{i,j}(D) \cong H^u_v(\mathfrak{C}(D),\partial),
\]
where both homologies are computed over $\zed$, $H^u_v(\mathfrak{C}(D),\partial)$ is the subspace of\break $H(\mathfrak{C}(D),\partial)$ of homogeneous elements of bidegree $(u,v)$, and
\begin{eqnarray*}
u & = & j-i-w(D)+1, \\
v & = & j-2i +\frac{E_+(G)+E_-(G)-w(D)}{2}+1.
\end{eqnarray*}
\end{theorem}

\begin{remark}
Since the $v$--grading in \cite{CKVST} is sensitive to the choice of coloring, the construction there is done under the assumption that $E_+(G)\geq E_-(G)$. But we are using a different normalization of the $v$--grading, which is invariant under reversing of the coloring. So \fullref{STmodel} is true for either coloring. The equivalence of \fullref{STmodel} and the corresponding result in \cite{CKVST}, and the invariance of our $v$--grading under reversing of the coloring, can be easily deduced from the following discussion of dual graphs.
\end{remark}

\begin{definition}\label{dualgraphs}
Let $G$ be a graph embedded in $\Real^2$. The dual graph $G'$ of $G$ is a graph embedded in $\Real^2$ defined as following:

(i)\qua All vertices of $G'$ are in $\Real^2\setminus G$. And each connected component of $\Real^2\setminus G$ contains exactly one vertex of $G'$.

(ii)\qua There is a one-to-one correspondence between edges of $G$ and $G'$, called the dual relation, under which any edge $e$ of $G$ corresponds to an edge $e'$ of $G'$ that transversally intersects $e$ once, connects the vertices of $G'$ in the connected components of $\Real^2\setminus G$ on both sides of $e$, and is disjoint from all other edges of $G$.

If the edges of $G$ are signed, then the edges of $G'$ are signed so that dual edges have opposite signs. If the edges of $G$ are ordered, then we order the edges of $G'$ so that the dual relation of edges preserves the ordering of edges.
\end{definition}

Note that the dual of the dual of a graph is the original graph. Also, the two Tait graphs of a link diagram obtained by reversing the coloring are duals of each other.

\begin{figure}[ht!]
\centering
\includegraphics[width=2in]{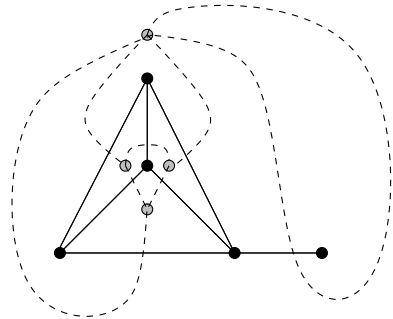}
\caption{Dual Graphs}
\end{figure}

\begin{lemma}\label{dualtrees}
Let $G$ be a graph embedded in $\Real^2$, and $G'$ its dual. Assume that both $G$ and $G'$ are connected. Then there is a one-to-one correspondence between spanning trees of $G$ and $G'$, which is called the dual relation, so that, for any spanning tree $T$ of $G$, its dual spanning tree $T'$ of $G'$ is defined by $e'\in T'$ if and only if $e\notin T$, where $e'$ is the dual of $e$.

Moreover, if edges of $G$ are signed and ordered, and the edges of $G'$ are signed and ordered as in \fullref{dualgraphs}, then, for any spanning tree $T$ of $G$ with dual spanning tree $T'$ of $G'$, and edge $e$ of $G$ with dual edge $e'$, we have
\begin{eqnarray*}
e \text{ is of type } L \text{ under } T & \Leftrightarrow & e' \text{ is of type } \overline{l} \text{ under } T'; \\
e \text{ is of type } D \text{ under } T & \Leftrightarrow & e' \text{ is of type } \overline{d} \text{ under } T'; \\
e \text{ is of type } l \text{ under } T & \Leftrightarrow & e' \text{ is of type } \overline{L} \text{ under } T'; \\
e \text{ is of type } d \text{ under } T & \Leftrightarrow & e' \text{ is of type } \overline{D} \text{ under } T'.
\end{eqnarray*}
\end{lemma}

\begin{proof}
We call a subgraph of $G$ a spanning subgraph if it contains all the vertices of $G$. For any spanning subgraph $H$ of $G$, define its dual spanning subgraph $H'$ of $G'$ by $e'\in H'$ if and only if $e\notin H$ for any pair of dual edges $e$ and $e'$. This is a one-to-one correspondence between spanning subgraphs of $G$ and $G'$, called the dual relation. We need to show that a spanning subgraph of $G$ is a tree if and only if its dual is a tree.

We compactify $\Real^2$ to $S^2$ by adding a single point at $\infty$. Let $H$ be a spanning subgraph of $G$, and $H'$ its dual spanning subgraph of $G'$. Slightly thicken $H$ in $S^2$. We get a surface $U_H$ with boundary which has $H$ as a deformation retract. Let $V_{H'}$ be the closure of $S^2\setminus U_H$ in $S^2$. Then $V_{H'}$ is a surface with boundary, and has $H'$ as a deformation retract. Note that $U_H \cup V_{H'}=S^2$, and $U_H \cap V_{H'}=\partial U_H = \partial V_{H'}$. It's clear that:

$H$ is a spanning tree. $\Leftrightarrow$ $U_H$ is a disk. $\Leftrightarrow$ $V_{H'}$ is a disk. $\Leftrightarrow$ $H'$ is a spanning tree.

Thus, the dual relation gives a one-to-one correspondence between spanning trees of $G$ and $G'$.

Now let $T$ be a spanning tree of $G$ and $T'$ its dual spanning tree of $G'$. Let $e$ be an edge of $T$, and $e'$ its dual edge. Then $e'\notin T'$. Let $f$ be an edge of $G$ not contained in $T$. Then $f'$ , the dual of $f$, is contained in $T'$. If $f\in \cut(T,e)$, then $e\in \cyc(T,f)$, ie, $e$ is contained in the unique simple cycle in $T\cup f$. So $f'$ is the unique edge in $T'$ connecting vertices on two sides this cycle. Thus, the two connected components of $T'\setminus f'$ are on each side of this cycle. Note that $e'$ connects vertices on two sides of this cycle. So $e'\in \cut(T',f')$, ie, $f'\in \cyc(T',e')$. If $f'\in \cyc(T',e')$, then $e'\in \cut(T',f')$, and, by the above argument, one can check that $e\in \cyc(T,f)$, ie, $f\in \cut(T,e)$. Thus, the dual relation of edges gives a one-to-one correspondence between $\cut(T,e)$ and $\cyc(T',e')$. This correspondence implies the second half of the lemma.
\end{proof}

\section{Legendrian links and Khovanov homology}\label{legendrian}

All Legendrian links in this paper are in the standard contact $\Real^3$, which is defined by the contact form $\alpha=dz-ydx$.

The front diagram of a Legendrian link is its projection onto the $xz$--plane. It's a immersion of circles into $xz$--plane with cusps and transversal self-intersections, but no vertical tangencies. Using the equation $y=\frac{dz}{dx}$, it's easy to check that a Legendrian link is uniquely determined by its front projection. After a small perturbation, we assume that all self-intersections are transversal double points (crossings) with pairwise different $x$--coordinates. In the rest of this paper, we will order the crossings of a front diagram by their $x$--coordinates so that the order of crossings increases from left to right. It's easy to see that, at any double point of the front diagram, the branch with smaller slope is on top.

\begin{figure}[ht!]
\centering
\includegraphics[width=.75in]{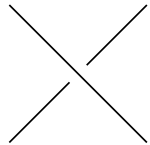}
\caption{The branch with smaller slope at a crossing is on top}\label{frontcrossing}
\end{figure}

Comparing Figures \ref{sign} and \ref{frontcrossing}, it's easy to see that a crossing is positive if and only if it is horizontal in the sense that each of the two black regions incident on the crossing is above one branch and below the other, and a crossing is negative if and only if it is vertical in the sense that one of the two black regions incident on the crossing is above both of the branches and the other is below both of the branches.

\begin{figure}[ht!]
\centering
\labellist\small
\pinlabel $+$ [b] at 181 42
\pinlabel $-$ [l] at 435 40
\endlabellist
\includegraphics[width=4in]{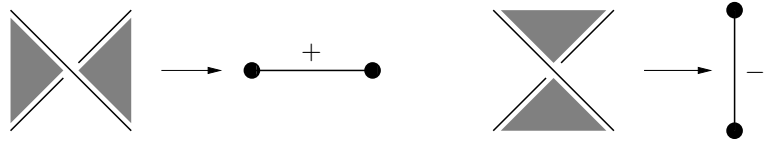}
\caption{Horizontal and vertical edges}\label{vh}
\end{figure}

\begin{figure}[ht!]
\centering
\labellist
\pinlabel $A$ [b] at 113 38
\pinlabel $B$ [b] at 251 38
\endlabellist
\includegraphics[width=3in]{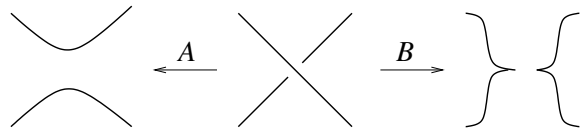}
\caption{Legendrian $A$-- and $B$--splicings}\label{legendriansplicings}
\end{figure}

We define the Legendrian $A$-- and $B$--splicings of a crossing in a front diagram by \fullref{legendriansplicings}.

Let $F$ be the front diagram of an oriented Legendrian link. Let $D$ be the desingularization of $F$, ie, the oriented link diagram obtained from $F$ by smoothing all the cusps. Then $F$ and $D$ represent the same topological link. They also have the same signed Tait graph $G$. Denote by $C(F)$ half of the number of cusps in $F$. The Thurston--Bennequin number of $F$ is $tb(F)=w(D)-C(F)$.

In the rest of this section, we assume that the diagram $F$ is connected as a subset of $\mathbb{R}^2$. (This can be easily arranged using Legendrian Reidemeister moves.)

For a spanning tree $T$ of $G$, let $U_T$ be the unknot obtained from
$D$ and $T$ using \fullref{splicingdead}, and $F_T$ the front
diagram obtained from $F$ and $T$ using \fullref{splicingdead} while
interpreting "$A$" and "$B$" as Legendrian $A$-- and
$B$--splicings. Note that $F_T$ is a Legendrian unknot whose
desingularization is $U_T$.

\begin{proposition}\label{tbFT}
If the front $F$ is connected, then, for any spanning tree $T$ of $G$,
\begin{equation}\label{FT}
tb(F_T) \leq -1-(\#d(T)+\#\overline{D}(T)),
\end{equation}
which is equivalent to
\begin{equation}\label{FTalter}
u(T) \geq 1- C(F).
\end{equation}
\end{proposition}
\begin{proof}
We prove Inequality \eqref{FT} by induction on the number of crossings in $F$. When $F$ has no crossings, $G$ consists of a single vertex. So the only spanning tree of $G$ is $T=G$, which satisfies $\#d=\#\overline{D}=0$. Note that $F_T=F$ is a Legendrian unknot. By Bennequin's inequality in \cite{Ben}, we have $tb(F_T) \leq -1 = -1-(\#d+\#\overline{D})$. This shows that the proposition is true for any front diagram with no crossings.

Now assume that \eqref{FT} is true for any connected front diagram with less than $m$ crossings. Let $F$ be a connected front diagram with $m$ crossings. Let $e$ be the edge of $G$ corresponding to the right most crossing of $F$. (In the rest of this proof, we do not distinguish between the edge $e$ and the crossing it represents.) $e$ has the highest order of all edges of $G$. There are three possibilities:
(1) $e$ is an isthmus, ie, $G$ becomes disconnected after removing $e$ from it; (2) $e$ is a loop, ie, $e$ connects a vertex to itself; (3) $e$ is neither an isthmus nor a loop.

\medskip
{\bf Case ($1$)\qua $e$ is an isthmus}\qua Then any spanning tree $T$ contains $e$, and $\cut(T,e)=\{e\}$. So $e$ is always internally active. Let $G_1$ and $G_2$ be the two connected components of $G\setminus e$. There are two possibilities:

\medskip
{\bf Case ($1_1$)\qua $e$ is negative and, hence, vertical}\qua The contribution of such an $e$ to the Thurston--Bennequin number is $+1$. We $A$--splice $F$ at $e$, which gives us two disjoint connected front diagrams $F_1$ and $F_2$, with Tait graphs $G_1$ and $G_2$. Any spanning tree $T$ of $G$ is the union of $e$ and a spanning tree $T_1$ of $G_1$ and a spanning tree $T_2$ of $G_2$. Let $F_{T_i}$ be the Legendrian unknot obtained from $F_i$ by the spanning tree $T_i$, $i=1,2$. Then, $tb(F_T) = tb(F_{T_1}) + tb(F_{T_2}) +1$, $\#d(T) = \#d(T_1) + \#d(T_2)$ and $\#\overline{D}(T) = \#\overline{D}(T_1) + \#\overline{D}(T_2)$. By induction hypothesis, we have $tb(F_{T_i}) \leq -1-(\#d(T_i)+\#\overline{D}(T_i))$, $i=1,2$, which implies that
$tb(F_T) \leq -1-(\#d(T)+\#\overline{D}(T))$.

\medskip
{\bf Case ($1_2$)\qua $e$ is positive and, hence, horizontal}\qua The contribution of such an $e$ to the Thurston--Bennequin number is $-1$. We $B$--splice $F$ at $e$, which gives us two disjoint connected front diagrams $F_1$ and $F_2$. Note that the sum of the contributions of the two new cusps to the Thurston--Bennequin number is also $-1$. Let $T$, $T_i$, $F_{T_i}$ be similarly defined as in ($1_1$). Then $tb(F_T) = tb(F_{T_1}) + tb(F_{T_2})$, $\#d(T) = \#d(T_1) + \#d(T_2)$ and $\#\overline{D}(T) = \#\overline{D}(T_1) + \#\overline{D}(T_2)$. By induction hypothesis, we have $tb(F_{T_i}) \leq -1-(\#d(T_i)+\#\overline{D}(T_i))$, $i=1,2$, which implies that $tb(F_T) \leq -2-(\#d(T)+\#\overline{D}(T)) < -1-(\#d(T)+\#\overline{D}(T))$.

\medskip
{\bf Case ($2$)\qua $e$ is a loop}\qua Let $G'$ be the Tait graph from the other coloring of $D$, which is the dual of $G$. Then the edge $e'$ of $G'$ dual to $e$ is an isthmus. By \fullref{dualtrees}, one can easily reduce this case to Case (1).

\medskip
{\bf Case ($3$)\qua $e$ is neither an isthmus nor a loop}\qua In this case, $e$ is inactive (internally or externally) for any spanning tree $T$ of $G$ since it has the highest order among all crossings. We need to consider the type of $e$ with respect to the spanning tree $T$.

\medskip
{\bf Case ($3_1$)\qua $e$ is of type $\overline{d}$}\qua Then $e$ is negative and, hence, vertical. We $A$--splice $F$ at $e$, which gives a connected front $\hat{F}$. Let $\hat{G}$ be the corresponding Tait graph of $\hat{F}$. Then $\hat{G}=G\setminus e$, and $T$ is a spanning tree of $\hat{G}$. Let $\hat{F}_T$ be the front obtained from $\hat{F}$ using $T$. Note that any edge $\hat{e}$ of $\hat{G}$ has the same type in $\hat{G}$ under $T$ as in $G$ under $T$. So $F_T = \hat{F}_T$, $\#d_G(T) = \#d_{\hat{G}}(T)$ and $\#\overline{D}_G(T) = \#\overline{D}_{\hat{G}}(T)$. By induction hypothesis, we have $tb(\hat{F}_T) \leq -1 - \#d_{\hat{G}}(T)+\#\overline{D}_{\hat{G}}(T)$, and, therefore, $tb(F_T) \leq -1-(\#d_G(T)+\#\overline{D}_G(T))$.

\medskip
{\bf Case ($3_1'$)\qua $e$ is of type $D$}\qua Consider other Tait graph of $F$, which is dual to $G$. By \fullref{dualtrees}, one can easily reduce this case to Case ($3_1$).

\medskip
{\bf Case ($3_2$)\qua $e$ is of type $d$}\qua Then $e$ is positive and, hence, horizontal. We $B$-splice $F$ at $e$, which creates a pair of cusps, one of which opens to the right, the other opens to the left. Let $\hat{F}$ be the resulting front diagram, which is connected. Denote by $c$ the new right opening cusp. Since $e$ is the right most crossing, which is neither an isthmus nor a loop, the two branches intersecting at $c$ do not intersect elsewhere. Next, we use an observation made by Ng in \cite{NgKB}: Extend the two branches at $c$ along $\hat{F}$ in both directions until it passes to the left of $c$. The result is a zigzag that does not intersect other parts of $\hat{F}$. Let $c_1$ and $c_2$ be the two consecutive cusps on this zigzag such that the difference of the $x$--coordinates of these two cusps is the smallest among all pairs of consecutive cusps on this zigzag. This minimality forces the part of the zigzag near these two cusps to look like one of the two zigzags in \fullref{zigzags}. Thus, we can destabilize $\hat{F}$ by smoothing out these two cusps (cf Figure 1 of \cite{NgKB}). Let $\tilde{F}$ be the result. Note that the Tait graph $\tilde{G}$ of $\tilde{F}$ is the same as that of $\hat{F}$, which is $G\setminus e$.

\begin{figure}[ht!]
\centering
\includegraphics[width=2in]{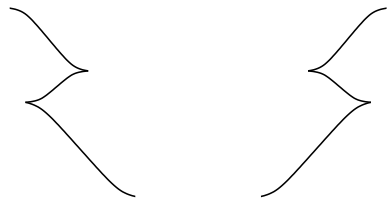}
\caption{The zigzag near $c_1$ and $c_2$}\label{zigzags}
\end{figure}

The spanning tree $T$ of $G$ is also a spanning tree of $\tilde{G}$. Let $\tilde{F}_T$ be the front obtained from $\tilde{F}$ using $\tilde{G}$ and $T$. It's easy to see that $tb(F_T) = tb(\tilde{F}_T)-1$, $\#d_G(T) = \#d_{\tilde{G}}(T)+1$ and $\#\overline{D}_G(T)= \#\overline{D}_{\tilde{G}}(T)$.
By the induction hypothesis, we have $tb(\tilde{F}_T) \leq -1-(\#d_{\tilde{G}}(T)+\#\overline{D}_{\tilde{G}}(T))$, and, therefore, $tb(F_T) \leq -1-(\#d_G(T)+\#\overline{D}_G(T))$.

($3_2'$) $e$ is of type $\overline{D}$. Consider other Tait graph of $F$, which is dual to $G$. By \fullref{dualtrees}, one can easily reduce this case to Case ($3_2$).

This completes the induction.

Note that
\begin{eqnarray*}
\#d(T)+\#\overline{D}(T) & = & \text{ number of Legendrian $B$--splicings used to obtain $F_T$} \\
& = & C(F_T)- C(F),
\end{eqnarray*}
and
\[
u(T) = -w(U_T) = -tb(F_T)-C(F_T).
\]
So Inequalities \eqref{FT} and \eqref{FTalter} are equivalent.
\end{proof}

\begin{definition}
Given a connected front diagram $F$ and an Tait graph $G$ of $F$, let $T$ be a spanning tree of $G$. $T$ is said to be good if
\[
u(T)=1-C(F),
\]
and is said to be bad if
\[
u(T)=2-C(F).
\]
For an integer $v$, $T$ is called a $v$--spanning tree if $v(T)=v$.
\end{definition}

\begin{proof}[Proof of Theorems \ref{NgBound} and \ref{sharpness}]
Let $F$ be a connected front projection of $L$, and $D$ the desingularization of $F$. Checkerboard color the front diagram, and let $G$ be the Tait graph from the coloring.

If $\mathcal{H}^{i,j}(L)\neq0$, then, by \fullref{STmodel}, there exists a spanning tree $T$ of $G$ such that $j-i=u(T)+w(D)\pm1$. By \fullref{tbFT}, we have $u(T) \geq 1- C(F)$. So
\[
j-i \geq w(D) - C(F)= tb(L).
\]
This proves \fullref{NgBound}.

If Ng's Khovanov bound is sharp for the Legendrian knot $L$, then there exists a pair $(i,j)$ such that $j-i=tb(L)$, $\mathcal{H}^{i,j}(L)\neq0$, and the corresponding spanning tree $T$ satisfies $j-i = u(T) + w(D)\pm 1 = tb(L)$, ie, $u(T)=\mp1-w(D)+tb(L)=\mp1-C(F)$. But, $u(T) \geq 1- C(F)$. So $u(T) = 1- C(F)$, which means $T$ is a good spanning tree.

Now assume $F$ has more good $v$--spanning trees than bad $(v+2)$--spanning trees for some integer $v$. Since the boundary map $\partial$ on $\mathfrak{C}(D)$ has bidegree $(-1,-2)$, this implies that $H^{1-C(F)}_v(\mathfrak{C}(D),\partial)$ has positive rank. So the corresponding pair $(i,j)$ satisfies $j-i=tb(L)$ and $\mathcal{H}^{i,j}(L)\neq0$. Hence, Ng's Khovanov bound is sharp for $L$.
\end{proof}

\begin{proof}[Proof of \fullref{alternating}]
Without loss of generality, we assume the alternating link $\mathcal{L}$ is non-split. In \cite{NgKB}, Ng constructed a Legendrian link $L$ with the link type $\mathcal{L}$ satisfying that, after choosing an appropriate checkerboard coloring for the front projection $F$ of $L$, all the crossings of $F$ are vertical (ie negative), and every black region is a bounded disk that has exactly two cusps on its boundary. Let $T$ be a minimal spanning tree of the Tait graph $G$ in the sense that the sum of the $x$--coordinates of the crossings corresponding to edges of $T$ is minimal among all the spanning trees of $G$. Then it is easy to check that all the edges of $T$ are internally active, and all the edges outside $T$ are externally inactive. Let $V(G)$ be the number of vertices of $G$. Then $V(G)=C(F)$, and
\begin{eqnarray*}
u(T) & = & \#L(T)-\#l(T)-\#\overline{L}(T)+\#\overline{l}(T) \\
     & = & -\#\overline{L}(T) = -(V(G)-1) = 1-C(F).
\end{eqnarray*}
So $T$ is a good spanning tree. Note that $F$ is an alternating diagram. From the proof of Theorem 12 of \cite{CKVST}, one can see that the $v$--grading of the spanning trees of $G$ is a constant $v_0$. Therefore, there are no $(v_0+2)$--spanning trees. Then \fullref{sharpness} implies that Ng's Khovanov bound is sharp for $L$.
\end{proof}

\bibliographystyle{gtart}
\bibliography{link}

\end{document}